\documentclass{article} 
\usepackage{graphicx} 
\begin{document}
\setcounter{page}{1} 

\begin{center}
{\Large\bf On the extreme complexity of certain nearly 
regular graphs } 
\end{center}
\medbreak
\centerline{by}
\medbreak
\centerline{Gregory P. Constantine}

\centerline{School of Computer Science}

\centerline{Georgia Institute of Technology}

\centerline{Atlanta, GA 30332 }

\medbreak
\centerline{and}
\medbreak
\centerline{Gregory C. Magda}

\centerline{Department of Mathematics}

\centerline{University of Pittsburgh}

\centerline{Pittsburgh, PA 15260}
\vskip1cm 
\centerline{\bf ABSTRACT}
 
\medbreak\noindent The complexity of a graph is the 
number of its labeled spanning trees.  It is demonstrated 
that the seven known triangle-free strongly regular 
graphs, such as the Higman-Sims graph, are graphs of 
maximal complexity among all graphs of the 
same order and degree; their complements are shown to 
be of minimal complexity.  A generalization to nearly 
regular graphs with two distinct eigevalues of the 
Laplacian is presented.  Conjectures and applications of these 
results to biological problems on neuronal 
activity are described.  

\noindent {\em Proposed }
{\em running head:\/}  Complexity of nearly regular graphs 
\vskip1cm 
\footnoterule\noindent Funded under NIH grant 
RO1-HL-076157 and NSF award ID 2424684 

\noindent gmc@pitt.edu 
  
\newpage

\section*{\Large\bf1.  Preliminaries} \medbreak\ The graphs 
in this work are finite, loopless, undirected, and without 
multiple edges.  They are also connected.  By the {\em order}
of a graph we understand the number of its vertices, 
and by {\em size} the number of its edges.  Standard 
terminology is used and we assume that the reader is 
familiar with such notions as path, graph connectivity, 
tree and spanning tree, adjacency matrix and the 
Laplacian.  For clarity we remind that a $walk$ of length 
$k$ (or $k-$walk) is a sequence of vertices and edges 
$v_1e_1v_2e_2\cdots v_ke_kv_{k+1}$, where $e_j$ is the edge joining vertices 
$v_j$ and $v_{j+1}.$ Vertices and edges may be repeated in this 
sequence.  The walk is {\em closed\/} if $v_1=v_{k+1}.$ An {\em m-cycle\/} is 
a sequence of vertices and edges $v_1e_1v_2e_2\cdots v_me_mv_{m+1}
,$ 
where all vertices $v_i$ are distinct except for $v_1$ and $v_{m+
1}$ 
which are the same; $m\geq 2.$ A $triangle$ is a $3-$cycle.  A 
graph is said to have {\em maximal complexity\/} if it has a 
maximal number of spanning trees among all graphs of 
the same order and size.  We aim, as far as we can, to 
understand the structure of such graphs. The results of 
this paper only achive this optimization among 
subclasses of graphs whose degrees differ by at most 1. Specific 
conjectures are then formulated and an asymptotic 
motivation is given to describe the general graph structure 
for graphs of 
maximal complexity among all graphs of a given order 
and size.
 
\vskip.3cm\noindent
Denote by $D$ the diagonal matrix with the degrees of the 
vertices of graph $G$ as entries (written always in the 
same fixed order), by $A$ the adjacency matrix and by 
$L=D-A$ the Laplacian.  Recall that $A^r$ has the number 
of walks of length $r$ that start and end at vertex $v_i$ as 
the $i^{th}$ diagonal entry.  Kirchhoff's theorem states that 
the complexity of a graph (up to sign) is equal to the 
determinant of any principal minor of $L$ obtained by 
deleting a row and column of $L.$ It is an immediate 
consequence of Kirchhoff's theorem that {\em the number of }
{\em spanning trees\/} (or the complexity) {\em of} $G${\em , which we denote by }
$t(G)${\em , is equal to} $\frac 1{(n+1)}x_1x_2\cdots x_n$ {\em where} $
(n+1)$ {\em is the }
$order$ {\em of} $G$ {\em and} $(0=)x_0\leq x_1\leq x_2\leq\ldots
\leq x_n$ {\em are the }
{\em eigenvalues of the Laplacian} $L$ {\em of} $G.$ 
\vskip.3cm\noindent
For simplicity of notation we write $e_n=x_1x_2\cdots x_n.$ 
\vskip.3cm\noindent
In spite of the fact that the complexity of a graph is 
expressed as the determinant of (any) principal minor of 
the Laplacian, this does not easily reveal the actual 
graph features that optimize complexity; embarassingly 
perhaps, it is still not known whether graphs of 
maximal complexity are found among the set of graphs 
whose degrees differ by at most one.  [This seems 
intuitive since any vertex appears to hold equal 
importance in building spanning trees.]  We present 
detailed specific structural conjectures later in the 
paper.  \medbreak\section*{\Large\bf2.  Applications of 
graph complexity} \medbreak\noindent The study of graphs 
of maximal complexity, apart from being an attactive 
purely graph theoretical incentive, has several practical 
motivations.  We briefly narrate several non-exhaustive 
instances below.  \medbreak\noindent{\bf 1.  Optimal credence ellipsoids} In statistical inference, for a random vector $X$ 
with density $f(x,\theta ),$ the Fisher information matrix is 
\medbreak\noindent
\centerline{$L=(f_{ij}(\theta ))=E[(\frac {\partial}{\partial\theta_
i}lnf(x,\theta ))(\frac {\partial}{\partial\theta_j}lnf(x,\theta 
))].\mbox{\rm \           (1)}$}
\medbreak\noindent Matrix $L$ is positive semidefinite and 
defines a metric on the parameter space indexed by $\theta ;$ it 
is the Hessian of entropy.  Its study is central in 
establishing efficient strategies of parameter estimation 
in both Bayesian and frequentist contexts.  The 
information matrix (1) may be viewed as a generalization of 
the graph Laplacian.  Efficient statistical estimation 
strategies entail maximization of certain types of 
functions defined on the spectrum of $L;$ one of these is 
the inverse of the product of eigenvalues, which is 
proportional to the volume of the confidence 
hyper-ellipsoids that arise in estimation. The smaller the volume of the hyper-ellipsoids the better the estimation. 
\medbreak\noindent{\bf 2.  A simple additive response model }
\medbreak\noindent Under Gaussian assumptions of 
independent and identically distributed responses, the 
normal equations of an additive model lead to Fisher 
information matrices as inverses of covariance matrices 
for the estimates of the parameters of interest.  The 
simplest case is that of observations made on graphs, 
where the expected value of an observation made on 
vertex $i$ that lies on edge $j$ is $\theta_i+\beta_j$.  The Fisher 
information matrix (1) associated to the parameters $\theta_i$ is in 
this case exactly the Laplacian of the graph in question.  
Typically we seek to find confidence (or credence) 
ellipsoids for orthogonal functionals of $\theta =(\theta_i)$ of 
smallest volume.  The volume of such a confidence 
ellipsoid is known to be inversely proportional to the 
number of labeled spanning trees of the graph; see [6].  
\medbreak\noindent{\bf 3.  Network connectivity }
\medbreak\noindent In electrical engineering it is desirable 
to construct a network that has a maximal number of 
spanning trees, since this provides the largest number of 
routes for transmitting information from any node to 
any other, thus minimizing the chance of overall failure.  
Building cost being typically a factor, the optimization 
problem becomes that of constructing a graph with a 
maximal number of spanning trees, given the number of 
nodes and how much we can afford to spend on the links 
joining these nodes.  \medbreak\noindent{\bf 4.  A purely graph-theoretic measure of connectivity }
\medbreak\noindent In graph theoretic terms, we can 
measure the connectivity of a graph by taking a 
spanning tree as a unit of connectivity.  For a given 
number of vertices and a given number of edges, a graph 
with a maximal number of spanning trees is then clearly 
the most connected.  As we shall shortly conjecture, a 
graph of maximal complexity is nearly regular, which 
means that the degrees of its vertices are either $d$ or 
$d+1.$ All evidence on graphs with 12 or fewer vertices, 
and verification in certain infinite families, points to the 
fact that such graphs are in fact $d-$edge connected 
(which is clearly best possible); in other words they 
appear to be most resistant in remaining connected upon 
removal of edges.  \medbreak\noindent{\bf 5.  Optimal spreading of information in neuronal modeling networks }
\medbreak\noindent
The motivating biological problem is to 
turn on all neurons in a brain, or part of a brain, by 
starting with a small subset of active neurons.  We 
view this activity as having {\em local component\/}s, which we 
want to turn on as fast as possible, and {\em global links }
between the local components which serve the purpose 
of efficiently integrating the local components in such a 
way that the entire brain becomes active as quickly as 
possible.  Further detailed information is found in $[10]$ 
and [9].  Intuitively we thus seek to determine those 
neuronal configurations, viewed as abstract networks, 
that spread the information most efficiently (fastest 
possible) to the whole brain.  We focus initially on 
modeling the local components and start by making some 
simplifying assumptions.  The basic working hypothesis 
is that a neuron is activated by receiving input from at 
least $t$ already active neurons connected to it.  Initially 
we make the assumption that the underlying graph that 
connects the neurons is regular.  Since we want a quick 
spread to activate the whole local area, it is intuitive 
that the best way of doing this is to avoid having short 
cycles, like triangles or squares, in the regular graph.  
If we have $n$ neurons, each of degree $d,$ the emerging 
optimization strategy is that we want to first restrict 
to having a minimal number of triangles, then among 
this subset of regular graphs to seek those that have a 
minimum number of closed walks of length 4 (like 
4-cycles), and proceed sequencially to closed walks of 
higher order.  A consequence of this paper is establishing 
a mathematical connection between the choice strategy we 
just described and nearly regular graphs that have 
a maximum number of spanning trees. The interplay 
between short cycles and spanning trees is highlighted 
later in the paper. We describe next, 
in some detail, measures of the spread of neuronal 
activity.  \medbreak\noindent Imagine for a moment that 
the vertices of the graph (or digraph) are neurons and 
any existing edge transmits information form one neuron 
to another.  We start with a set $S$ of neurons, which 
we call $active,$ and a startup treshold $t,$ which is a 
natural number.  The spreading of neuronal activity is 
described next.  This is subject to some restrictions 
formulated in terms of {\em Steps}, which we now describe.  
\vskip.3cm\noindent
$Step$ 0:  Start with a set $S=S_0$ of vertices of the 
digraph $G$ and a natural number $t.$ We call elements of 
$S$ active vertices.  [Imagine that you hold the active 
vertices in your left hand, and the other vertices in 
your right hand.]  Color any edge emanating from $S$ red.  
\vskip.3cm\noindent
$Step$ 1:  Acquire vertex $v$, held in your right hand, if $v$ 
has $t$ or more red arrows pointing to it.  Move all 
acquired vertices to your left hand.  Call the set of 
vertices you now hold in your left hand $S_1.$ Color all 
edges emanating from $S_1$ red.  
\vskip.3cm\noindent
The general step is as follows.  We are in posession of 
$S_{i-1}$ with all edges emanating from it colored red.  
\vskip.3cm\noindent
$Step$ $i:$ Acquire vertex $v,$ held in your right hand, if it 
has $t$ or more red arrows pointing to it.  Move all 
acquired vertices to your left hand.  Call the set of 
vertices you now hold in your left hand $S_i.$ Color all 
edges emanating from $S_i$ red.  
\vskip.3cm\noindent
Evidently $S=S_0\subseteq S_1\subseteq\cdots\subseteq S_i\subseteq
\cdots .$ As we keep 
increasing $i$, the following will (obviously) always occur:  
the number of vertices in your right hand becomes 
stationary; that is, $\exists m$ such that, at Step $i,$ for all $
i\geq m$ 
the number of vertices in your right hand remains 
constant.  If your right hand becomes empty for a 
sufficiently large $i$ we say that the network is in 
{\em synchrony}.  [You are now holding the whole network in 
your left hand -- hence all vertices of the network 
became active.]  We denote by $i^{*}[=i^{*}(S,t)]$ the smallest 
$i$ such that at Step $i$ the network is in synchrony.  
Typically a network cannot be brought to synchrony 
(starting with an incipient set $S=S_0$ and $t)$, and we 
convene to write $i^{*}=i^{*}(S,t)=\infty$ in such a case.  
\medbreak\noindent Fix $t.$ In a graph (or digraph) $G,$ let 
$S=S_0$ be a set with $k$ vertices, which we call a 
$k-$subset.  Write $i^{*}(S)$ for $i^{*}(S,t).$ We introduce the 
following measures of synchrony for $G$ and $k.$ 
\vskip.3cm\noindent
The ratio $p_k(G)=$(number of $k-$subsets $S$ that bring $G$ to 
synchrony)/(number of all $k-$subsets) signifies the 
probability of bringing digraph $G$ to synchrony from a 
randomly chosen $k-$subset.  Generally we are interested 
in identifying digraphs with large $p_k$.  It might also be 
observed that there are many instances when a digraph 
has a large $p_k$ but the number of steps required to 
obtain synchrony are generally quite large, which is not 
so good.  We could tune this up by defining another 
measure $e_k(G)$, which we call $synchrony$ $efficiency,$ as 
follows:  
\vskip.3cm 
${n\choose k}e_k(G)=\sum_S(i^{*}(S))^{-1}.\mbox{\rm \                         (2)}$
\vskip.3cm\noindent
Observe that when $S$ does not induce synchrony, 
$i^{*}(S)=\infty ,$ and we simply add a zero to the sum.  
Intuitively, efficiency $e_k$ yields the average speed to the 
synchrony of $G$ across all $k-$subsets.  High values of $e_k$ 
are typically good, since the synchrony is then speedily 
restored.  We did not see the concept of synchrony 
efficiency, or analogs of (2), used in the network optimization literature so 
far. 
\medbreak\noindent
{\bf 6.  Tree expanders }
\medbreak\noindent
Recall that we denote by $t(G)$ 
the number of labeled spanning trees in the graph $G.$ 
For a sequence of graphs $G_n$ we write $t_n:=t(G_n).$ 
Furthermore, when we fix a degree $d$ and look at a 
sequence of graphs $(G_n)$, with each $G_n$ regular of degree $d$ 
for all $n$, we write $t_{d,n}$ for $t(G_n)$ -- a quantity 
that typically depends on both $d$ and $n.$
\medbreak\noindent$Definition$ A tree-expander of degree $d$ 
is a sequence $(G_n)$ of graphs, each of degree $d$, such 
that $liminf_{n\rightarrow\infty}\sqrt[{n-1}]{t_{d,n}}=:c_d$ and $
c_d\geq 1$.
\medbreak\noindent
As is the case with general (spectral) expanders, $c_d$ intuitively measures the 
connectivity at infinity of the sequence of $d$-regular 
graphs $(G_n)$. More specifically, the $c_d$ defined above 
measures the connectivity at infinity by {\em the geometric mean of the number of spanning trees\/} carried to 
infinity by the tree-expander $(G_n).$ [In particular we 
insist on $c_d\geq 1$ since we want on (geometric) average 
at least one spanning tree at infinity.] Evidently a ''good'' tree-expander for a given degree $d$ is 
a sequence $(G_n)$ for which $c_d$ is as high as possible; this 
happens, for example, when $G_n$ is a graph of maximal complexity, $
\forall n.$ 
Such graphs are generally speaking quite difficult to characterize, but sequences 
that consist of graphs of relatively high complexity may 
be constructed. To set the tone of inquiry, observe that for $d=2$ the 
sequence $(C_n)$, where $C_n$ is the cycle with $n$ vertices, yields 
$lim_{n\rightarrow\infty}x_1(C_n)=0<lim_{n\rightarrow\infty}\sqrt[{
n-1}]{t(C_n)}=lim_{n\rightarrow\infty}\sqrt[{n-1}]{n}=1$. 
This shows that $(C_n)$ {\em is not\/} a (usual) spectral-expander but $is$ a 
tree-expander. [We elucidate that it is easy to verify that 
$x_1(C_n)=2-2cos(\frac {2\pi}n)$ and $t(C_n)=n$, which explains both 
limits.]
\medbreak\noindent
Assume now that $d=3.$ We make tree-expanders the graphs $(M_{2n}
)$, where 
$M_{2n}$ is the M\"obius ladder with $2n$ vertices; $2n\geq 6.$ 
It is straightforward enough to show 
that $t(M_{2n})=\frac n2[(2+\sqrt {3})^n+(2-\sqrt {3})^n+2]\rightarrow\sqrt {
2+\sqrt 3}=c_3$ as $n\rightarrow\infty .$ 
\medbreak\noindent
More generally now, assume that we have a tree-expander $(H_n)$ for degree $
d.$ We 
shall construct a tree-expander $(G_n)$ of degree $d+1.$ 
Write $V_n$ for the number of vertices in $H_n$ and $E_n$ for 
the number of edges in $H_n.$  
Start with two copies on $H_n;$ call them $H_n^1$ and $H_n^2.$ 
Produce new edges $\{v_i^1,$ $v_i^2$$\}$ by joining vertex vertex $
v_i^1$ in 
$H_n^1$ to the corresponding vertex $v_i^2$ in $H_n^2$, $\forall 
i.$ The new 
graph, which we call $G_n$, has as vertices the union of 
vertices in $H_n^1$ and $H_n^2$ and as edges the edges in 
$H_n^1$, the edges in $H_n^2$, plus the newly created edges $\{v_
i^1,$ $v_i^2$$\}$, $\forall i.$
The number of vertices in $G_n$ is $2V_n$; the number of edges 
in $G_n$ is $2E_n+V_n.$ It can be shown that $(G_n)$ is a 
tree-expander and that $c_{d+1}>c_d.$ Details are omitted and 
shall appear elsewhere. 
\medbreak\section*{\Large\bf3.  
Optimizing complexity over nearly-regular graphs} 
\medbreak\medbreak\noindent For all $r=0$, define as 
follows:  \medbreak
\centerline{${\it g}_0=$ $\{$all connected graphs on $n+1$ vertices having at most $
m$ edges$\}$, }
\medbreak\noindent
and, for $r\geq 1$,
\medbreak
\centerline{${\it g}_r=\{G\in {\it g}_{r-1}:(-1)^{r-1}Tr(L^r)$ is maximal$
\}$. }
\vskip.5cm\noindent Evidently
\vskip.3cm
\centerline{ $_{g_0\supseteq {\it g}_1\supseteq {\it g}_2\supseteq
\cdots\supseteq {\it g}_r\supseteq\cdots}$ }
\vskip.3cm\noindent
and ${\it g}_r$ is {\em nonempty\/} for all $r\geq 0.$
\medbreak\noindent
Establish explicit constraint intervals $[a_r,b_r]$ by setting \hfill\break
\centerline{$a_r=min_{L\in {\it g}_{r-1}}Tr(L^r)$ and $b_r=max_{
L\in {\it g}_{r-1}}Tr(L^r)$, for $r\geq 1.$}
Observe, in particular, that $a_1=0,$ and 
$b_1=max_{L\in {\it g}_0}Tr(L)=2m,$ where $m$ is the number of edges 
of the graph whose Laplacian is $L$. 
\medbreak\noindent
In words, $g_0$ consists of graphs with $n+1$ vertices and at 
most $m$ edges. Set $g_1$ consists of graphs with $n+1$ 
vertices and exactly $m$ edges. Set $g_2$ consists of graphs 
with $n+1$ vertices that are nearly regular with each 
vertex of degree either $d$ or $d+1$, where $d$ is the 
integral part of $\frac {2m}{n+1}$. As mentioned above $g_0\supseteq 
g_1\supseteq g_2.$ 
Our approach is to optimize complexity by imposing 
restrictions on the trace functions of powers of 
Laplacians. We find such constraints quite intuitive, but 
it may be possible to use other bases of symmetric 
functions for this purpose; see [1].
\medbreak\noindent
For fixed $n$ and $d$ our aim is to identify graphs of maximal 
complexity in $g_2$. This is equivalent to maximizing $e_n$ over $
g_2.$
\medbreak\noindent 
Observe first that, by the monotonicity of the logarithm, 
$e_n$ is maximal if and only if $lne_n=\sum_{i=1}^nlnx_i$ is maximal.  
At this point we prefer to work with the function 
$f=f(x_1,\ldots ,x_n)=\sum_{i=1}^nlnx_i=lne_n.$ 
\medbreak\noindent
Intervals $[a_r,b_r]$ are as defined above. Being in ${\it g}_1$ assures us that the graphs have 
$n+1$ vertices and that $Tr(L)=\sum x_i=2m.$ We maximize $f$ 
subject to constraint $\sum x_i=2m(=b_1).$ By Lagrange 
multipliers the answer is that $f$ is maximal if and only 
if all $x_i$ are equal to $\frac {2m}n$.  Key issue to always raise is 
this:  
\vskip.3cm 
\centerline{{\em Are there any graphs with this Laplacian eigenstructure\/}?  }
\vskip.3cm\noindent
For arbitrary $m$ it is easy to see that this approach 
does not yield much information.  Specifically, a graph $G$ 
with this eigenstructure (that is with all $x_i$ equal to 
$\frac {2m}n$, $1\leq i\leq n$) exists if and only if $G=K_{n+1},$ the 
complete graph with $n+1$ vertices.  In particular this 
forces the number of edges $m$ to be $m={{n+1}\choose 2}$, and the 
degree $ $$d$ to be $n.$ [To see this, first note that $L$ 
satisfies its minimum polynomial $L^2-xL=0,$ with 
$x=\frac {2m}n$.  Look at the $i^{th}$ diagonal element of $L^2-\frac {
2m}nL=0$ 
to conclude that the graph is in fact regular of degree 
$d=\frac {2m}n-1$.  By regularity we must also have 
$2m=(n+1)d.$ Now $nd+d=2m=nd+n$ yields $d=n,$ 
forcing the graph to be complete.]  \medbreak\noindent The 
overarching conclusion is that the discretness of the set 
of graphs severely limits the spectral approach, at least 
initially.  The continuous constraint $\sum x_i=2m$ forces the 
''ideal'' eigenstructure with all $x_i$ equal as sole solution, 
and only the graph $K_{n+1}$ conforms to that.  
\vskip.5cm\noindent 
With that out of the way, let us focus on 
${\it g}_2.$ By definition, $a_2=min_{L\in {\it g}_1}Tr(L^2)$ is the minimal value 
that $Tr(L^2)$ takes among all graphs with $n+1$ vertices 
and $m$ edges.  It follows that ${\it g}_2$ consists of all graphs 
with $n+1$ vertices and $m$ edges whose degrees differ by 
at most 1.  These are the nearly regular graphs. [Evidently, regular 
graphs are nearly regular.]  We thus seek to maximize 
$f=\sum lnx_i$ subject to constrains $TrL=\sum x_i=2m=b_1$ and 
$Tr(L^2)=\sum x_i^2=a_2.$ Lagrange tells us that 
\[gradf=(\frac 1{x_1},\cdots ,\frac 1{x_n})=\lambda_1(1,\cdots ,1
)+\lambda_2(x_1,\cdots ,x_n).\mbox{\rm ~~~~~~~~~~(3)}\mbox{\rm }\]
This forces $\frac 1{x_i}=\lambda_1+\lambda_2x_i$, for all $i$.  The $
x_i$ are 
therefore roots of the quadratic $\lambda_2x^2+\lambda_1x-1=0$.  Note 
that $\lambda_2$ cannot be zero, else the $x_i$ would emerge as all 
equal again, and we would be in the previous case.  Thus 
the eigenvalues $x_i$ $(i\neq 0)$ of the Laplacian $L$ of $G$ take {\em exactly }
{\em two distict values.  }
\medbreak\noindent
For convenience we write $r=b_1=2m$ and $s=a_2.$ As 
mentioned, $s$ is the value of the sum of squares of 
degrees attained by any nearly regular graph with $n+1$ 
vertices and $m$ edges; such graphs always exist. Let $x_1$ 
and $x_2$ be the two distinct nonzero eigenvalues of $L$ with 
respective multiplicities $n_1$ and $n_2.$ We set $x_1>x_2>0.$ 
The constraints $n_1x_1+(n-n_1)x_2=r$ and $n_1x_1^2+(n-n_1)x_2^2=
s$ 
yield, with $z=\sqrt {ns-r^2}$,
\vskip.3cm
\centerline{$x_1=\frac rn+\frac zn\sqrt {\frac {n-n_1}{n_1}}$ and $
x_2=\frac rn-\frac zn\sqrt {\frac {n_1}{n-n_1}}$.       (4)}
\vskip.3cm\noindent
For fixed $(n,r,s)$, and hence fixed $z,$ it is not clear how many graphs with
such eigenstructure exist, if any. But the above 
expression for $x_1$ and $x_2$ indicates that they can be 
indexed by $n_1;$ for notational ease write $x$ for $n_1.$ 
We may now write our optimizing function as 
$f=xlnx_1+(n-x)lnx_2,$ and view it solely as a function of 
$x.$ Specifically,
\vskip.3cm
\centerline{$f=xln(r+\frac {(n-x)z}{\sqrt {x(n-x)}})+(n-x)ln(r-\frac {
xz}{\sqrt {x(n-x)}})-nln(n),$ $1\leq x\leq n-1.$}
\vskip.3cm\noindent
It can be shown, and the details are found in the 
Appendix, that \medbreak\noindent {\bf Lemma 1} {\em Function} $f$ {\em is strictly}
{\em decreasing in} $x,$ $1\leq x\leq n-1.$ \medbreak\noindent This tells 
us that of two nearly regular graphs, with the same 
order and size, and each having two distinct non-zero eigenvalues 
of the Laplacian, the one that has a lesser multiplicity 
of its larger eigenvalue has more spanning trees. To shed 
further clarity on how to identify the graphs of maximal 
complexity over ${\it g}_2$, we shall examine $Tr(L^3)=\sum x_i^
3$. It is 
easy to see that $\sum x_i^3=Tr(L^3)=\sum d_i^3+3\sum d_i^2-6\Delta 
,$ with 
$d_i$ being the degrees and $\Delta$ the number of triangles in 
the graph. Furthermore, since the graphs are in ${\it g}_2,$ 
$\sum d_i^3$ takes the same value (say $M_3)$, and $\sum d_i^2$ takes the 
same value $(M_2)$ for all such graphs. It follows that {\em for }
{\em any two graphs} $G_1$ {\em and} $G_2$ {\em in} ${\it g}_2$ {\em with Laplacians} $
L_1$ {\em and }
$L_2,$ {\em we have} $Tr(L_1^3)\geq Tr(L_2^3)$ {\em if and only if} $
\Delta_1\leq\Delta_2.$
\medbreak\noindent
For a graph in ${\it g}_2$ with Laplacian $L$ and eigenvalues 
$x_1>x_2$ with multiplicities $n_1=x$ and $n_2=n-x,$ define 
now the function 
\vskip.3cm
\centerline{$c=Tr(L^3)=xx_1^3+(n-x)x_2^3=x(r+\frac {(n-x)z}{\sqrt {
x(n-x)}})^3+(n-x)(r-\frac {xz}{\sqrt {x(n-x)}})^3$ }
\vskip.3cm\noindent
We proved, with details available in the Appendix, that 
\medbreak\noindent {\bf Lemma 2} {\em Function} $c$ {\em is strictly }
{\em decreasing in} $x,$ $1\leq x\leq n-1.$ \medbreak\noindent Using 
Lemmas 1 and 2 and employing a composition of 
bijections, along with the observation we made that 
$Tr(L^3)$ is decreasing in the number of triangles, we may 
now conclude as follows \medbreak\noindent {\bf Lemma 3} {\em On }
{\em graphs in} ${\it g}_2$ {\em with two distinct eigenvalues of the }
{\em Laplacian, function} $f$ {\em is strictly decreasing as a }
{\em function of the number of triangles in the graph.  }
\medbreak\noindent We now state and prove a result that 
idetifies graphs of maximal complexity over ${\it g}_2$.  Recall 
that ${\it g}_2$ consists of all graphs with $n+1$ vertices and $
m$ 
edges whose degrees differ by at most 1.  
\medbreak\noindent {\bf Theorem 1} {\em If a graph} $G^{*}$ {\em in} $
{\it g}_2$ {\em has }
{\em two distinct nonzero eigenvalues of the Laplacian and a }
{\em minimal number of triangles among all graphs in} ${\it g}_
2${\em , }
{\em than the graph} $G^{*}$ {\em has a maximal number of spanning }
{\em trees among all graphs in} ${\it g}_2.$ {\em Such a graph is unique }
{\em up to} $cospectrality${\em .  \medbreak}\noindent To prove this 
theorem, let $G^{*}$ be a graph that satisfies the enunciated 
assumtions in the theorem.  If $G\in {\it g}_2$ then $Tr(L)=2m$ and 
$Tr(L^2)=s,$ where $s$ is the sum of squares of the degrees 
of $G$ (which is nearly equal); note that $s$ is the same 
for all graphs in ${\it g}_2.$ Insisting that $G$ be of maximal 
complexity over ${\it g}_2$ forces $L$ to have exactly two distinct 
nonzero eigenvalues (by using the Lagrange multipliers 
argument).  For arbitrary $n$ and $m$ typically there are no 
graphs like that; but, by our assumption, we have $G^{*}$ in 
hand, which $is$ like that.  Our setting therefore allows 
such graphs to exist and, in accordance with Lagrange, 
we need only examine these to pick the optimal one 
among them.  The tool employed to that purpose is 
Lemma 3 which tells us that, from among them, we 
should select a graph with a minimal number of triangle.  
Graph $G^{*}$ is like that, and is thence of maximal 
complexity.  Uniqueness of $G^{*}$ also follows from Lemma 
3 since picking the smallest number of triangles 
identifies $n_1$ uniquely.  Knowing $n_1$ determines $n-n_1$ and 
hence the eigenvalues $x_1$ and $x_2$ of the Laplacian 
uniquely.  We conclude that $G^{*}$ is unique up to 
cospectrality.  This ends the proof.  \medbreak\noindent
Some consequences of Theorem 1 are discussed below.  
\vskip.5cm\noindent
{\bf 1.}  It is known that a regular graph having two nonzero 
eigenvalues of the Laplacian is necessarily strongly 
regular [2].  Furthermore, a strongly regular graph is 
determined, up to isomorphism, by its two eigenvalues 
and their multiplicities.  By Theorem 1, this allows us 
to conclude that if ${\it g}_2$ contains a strongly regular graph 
$G^{*}$ that has a minimal number of triangles over ${\it g}_2$, then 
$G^{*}$ is the graph of maximal complexity over ${\it g}_2$ and is 
unique up to isomorphism.  
\vskip.3cm\noindent
Let us turn now to the setting in which ${\it g}_2$ does not 
contain regular graphs.  A graph in ${\it g}_2$ has $n+1$ vertices 
and $m$ edges; denote by $d$ and $d+1$ the two distinct 
degrees.  Let further $x_1$ 
and $x_2$ be the two nonzero eigenvalues of the Laplacian $L$ of 
$G^{*}.$ Denote by $J$ the matrix with all entries 1.  Then $L$ 
satisfies the functional equation $(L-x_1I)(L-x_2I)=\frac {x_1x_2}{
n+1}J$, 
obtained from its minimum polynomial; see also [3, page 
33].  Write $x_1+x_2=p$ and $x_1x_2=q.$ Identifying a 
diagonal element on both sides of this functional equation 
we conclude that both $d$ and $d+1$ are roots of the 
quadratic $z^2-(p-1)z+\frac {nq}{n+1}=0.$ View this as a system in 
$p$ and $q.$ It yields $\frac q{n+1}=\frac {d(d+1)}n.$ Select any vertices $
i$ and 
$j$ in $G^{*}$ that are not joined by an edge.  Examining the 
$(i,j)^{th}$ entry of the bove functional equation yields that 
the number of vertices in $G^{*}$ joined to both $i$ and $j$ is 
$\frac q{n+1},$ and integer, necessarily.  Integrality of $\frac 
q{n+1}=\frac {d(d+1)}n$ 
forces $n$ to divide $d(d+1)$.  
\vskip.3cm\noindent
We shall not attempt to classify here all nearly regular 
graphs that satisfy this arithmetic condition.  It suffices 
for now to just observe that a graph $G^{*}$ that is a 
complete graph with some non-intersecting (shall we 
say, parallel) edges removed satisfies this condition.  We 
summarize:  \medbreak\noindent {\bf Corollary 1 }
\vskip.3cm\noindent
$(a)$ {\em If a strongly regular graph} $G^{*}$ {\em with a minimal }
{\em number of triangles over} ${\it g}_2$ {\em exists in} ${\it g}_
2${\em , then} $G^{*}$ {\em has }
{\em maximal complexity over} ${\it g}_2$, {\em and} $it$ {\em is unique up to }
{\em isomorphism.  }
\vskip.3cm\noindent
$(b)$ {\em If} $there$ {\em exists non-regular} $G^{*}\in g_2$ {\em of order} $
n+1$ {\em whose degrees }
{\em are} $d$ {\em or} $d+1${\em , and} $G^{*}$ {\em has two distinct nonzero eigenvalues }
{\em of the Laplacian then} $n$ {\em necessarily divides} $d(d+1)
.$ 
{\em Furthermore, if such} $G^{*}$ {\em also has a minimum number of }
{\em triangles, then} $G^{*}$ {\em is a graph of maximal complexity. }
{\em Nonexhaustive examples of such graphs are }
{\em complete graphs with some parallel edges removed.  }
\vskip.3cm\noindent It is not hard to verify that a 
complete graph with some parallel edges removed, as in 
Corollary 1 part (b), is of maximal complexity among all 
graphs of the same order and size.  \medbreak\noindent\ {\bf 2.  }
We highlight a notable special case of Corollary 1, part 
(a).  If ${\it g}_2$ contains triangle-free strongly regular graphs, then we 
conclude the following:  
\vskip.3cm\noindent
{\bf Corollary 2} {\em Each of the seven triangle-free strongly }
{\em regular graphs is a graph of maximal complexity over }
{\em the set of all regular graphs with the same number of }
{\em vertices and the same degree as itself.  These graphs }
{\em are:  the 5-cycle, the Petersen, Clebsch, }
{\em Hoffman-Sigleton, Gewirtz, Mesner and Higman-Sims }
{\em graphs.  }
\vskip.3cm\noindent
The graphs in Corollary 2 have quite large symmetry 
groups. Notably, the index-two subgroup of the 
symmetry group of the Higman-Sims graph is known as 
the Higman-Sims group; it is one of the 26 sporadic 
simple groups discovered by Higman and Sims in 1968; cf. 
[11].
\vskip.5cm\noindent
{\bf 3.}  Allow now for the possibility that the complete 
multipartite graph $G^{*}=K_{q,\ldots ,q}$, having $p$ parts of 
cardinality $q$ each, exists in ${\it g}_2$.  The Laplacian of $
G^{*}$ has 
two nonzero eigenvalues $x_1=pq$ of multiplicity $p(q-1)$ 
and $x_2=pq-q$ of multiplicity $p-1$.  This graph is also 
known to have a minimal number of triangles among all 
graphs in ${\it g}_2.$ By Theorem 1 it follows that 
\vskip.3cm\noindent
$Graph$ $K_{q,\ldots ,q}$ {\em has maximal complexity over the set of }
{\em graphs of the same order and same degree as itself.  }
\vskip.3cm\noindent
This yields an explicit infinite family of graphs of 
maximal complexity over ${\it g}_2$.  \medbreak\noindent{\bf 4.}  A 
simple example, given next, helps us understand why 
insistence upon having a minimal number of triangles is 
essential in yielding maximal complexity.  We look at 
graphs on 9 vertices with 18 edges.  The first contender 
is the Lattice graph $H$ on 9 vertices of degree 4, a 
strongly regular graph with eigenvalues of the Laplacian 
3 and 6, each of multiplicity 4.  Graph $H$ has 11664 
spanning trees.  It is very pleasing to the eye, in the 
sense that its symmetry group is of order at least 72.  
We do not believe any other graph in its class to be 
more symmetric.  To top it off, $H$ is also 
self-complementary!  Hands down, on looks alone, it 
should be the winner on complexity as well.  Graph $H$ is 
not of maximal complexity, however.  Consider the 
second contender, a graph $G$ with 9 vertices, regular of 
degree 4, with edges 15 16 17 18 25 26 27 28 35 37 38 39 
46 47 48 49 59 69.  $G$ $is$ the graph of maximal 
complexity.  With only 4 symmetries, it is less pleasing 
to the eye than $H,$ but its number of spanning trees is 
12480, exceeding that of $H.$ Surprisingly perhaps, graph $G$ 
has six distinct nonzero eigenvalues.  The key point is 
that $G$ has only two triangles, while $H$ has six.  
\medbreak\noindent We observe that the graph feature of 
having a minimum number of triangles is {\em essential\/} for 
graphs of maximal complexity over $g_2.$ In contrast, the 
fact that we insist on having two distinct non-zero 
eigenvalues of the Laplacian in Theorem 1 is actually 
non-essential and is merely used as a simplifying 
technical assumption to establish maximal complexity in 
certain settings.  More specific structural conjectures are found in 
the last section.  \medbreak\noindent We observe in 
passing that, in general, cospectral Laplacians 
corresponding to two graphs of maximal complexity, do 
not necessarily force the graphs to be isomorphic.  We 
exhibit two nonisomorphic cospectral graphs of {\em maximal }
{\em complexity\/} on 8 vertices and 18 edges.  Take $K_{4,4}$ with 
one part having vertices 1, 2, 3, 4.  Obtain $G$ by adding 
edges 12 and 34 to $K_{4,4}$ and $H$ by adding 12 and 56 to 
$K_{4,4}.$ Both Laplacians have spectrum 8, 6, 6, 4, 4, 4, 4.  
Each graph has complexity 9216, which is maximal.  They 
are nonisomorphic since the off-diagonals in $A^2$ differ by 
3 in $G$ and by 4 in $H$.  
\vskip.5cm\noindent
Let us turn our attention now to connected graphs of 
{\em minimal\/} complexity over the set ${\it g}_2$.  Particularly helpful 
is the operation of graph complementation; we write $\bar {G}$ 
for the complement of graph $G.$ Assume that both $^{}G$ and 
$\bar {G}$ are connected graphs with $n+1$ vertices; $G$ has $m$ 
edges.  Write $L$ and $\bar {L}$ for the Laplacians of $G$ and $\bar {
G}.$ 
First, a few helpful observations:  
\vskip.5cm\noindent\ 
$(a)$ $x_i$ {\em is a nonzero eigenvalue of} $L$ {\em if and only if }
$n+1-x_i$ {\em is a nonzero eigenvalue of} $\bar {L}.$ 
\vskip.3cm\noindent
This is easy to see, since $L+\bar {L}$ equals the Laplacian of 
the complete graph, and $x_i$ and $n+1-x_i$ share the same 
eigenvector.  
\vskip.3cm\noindent
$(b)$ $G$ {\em has two nonzero eigenvalues of the Laplacian if }
{\em and only if} $\bar {G}$ {\em does.  }
\vskip.3cm\noindent
This follows immediately from $(a).$ 
\vskip.3cm\noindent
$(c)$ $G$ {\em is nearly regular if and only if} $\bar {G}$ {\em is.  }
\vskip.3cm\noindent
$(d)$ $\Delta (G)+\Delta (\bar {G})=$ ${{n+1}\choose 3}-nm+\frac 
12\sum_i(d_i(G))^2$ 
\vskip.3cm\noindent
Here $\Delta (H)$ denotes the number of triangles in graph $H,$ 
and $d_i(H)$ the degree of vertex $i$ in $H.$ Assertion $(d)$ 
surfaces upon coloring the edges of $G$ red, those of $\bar {G}$ 
blue, and counting the number of monochromatic 
triangles in the resulting (bicolored) complete graph; see 
also [4].  
\vskip.3cm\noindent
$(e)$ $G$ {\em has a minimal number of triangles over the set of }
{\em nearly regular graphs to which} $G$ {\em belongs if and only }
{\em if} $\bar {G}$ {\em has a maximal number of triangles over the set of }
{\em nearly regular graphs to which} $\bar {G}$ {\em belongs.  }
\vskip.3cm\noindent
Assertion $(e)$ follows immediately from $(d).$ 
\medbreak\noindent A minimization of complexity using 
Lagrange multipliers, as written in (3), yields the 
information that such a graph should also ideally have 
just two nonzero eigenvalues of the Laplacian.  Lemmas 
1, 2, and 3 are also directly applicable.  We formulate 
our final statement by viewing Corollary 1 through the 
prism of complementation.  Complements of graphs in 
Corollary 1 part (b) yield disconnected graphs.  Finally, 
using assertions $(a)-(e)$ we may state the following 
\medbreak\noindent {\bf Corollary 3} {\em If a connected strongly }
{\em regular graph} $G^{*}$ {\em with a maximal number of triangles }
{\em over} ${\it g}_2$ {\em exists in} ${\it g}_2${\em , then} $
G^{*}$ {\em has minimal complexity }
{\em over} ${\it g}_2$, {\em and} $it$ {\em is unique up to isomorphism.  }
\medbreak\noindent In particular, the complement of each 
of the seven triangle-free strongly regular graphs 
mentioned in Corollary 2 is a graph of minimal 
complexity over the set of all graphs with the same 
number of vertices and the same degree as itself.  Along 
with the infinite families of the complements of the 
graphs described in Corollary 1, Corollary 3 answers, 
in part, a question raised by Peter Sarnak and Noga 
Alon; see [5].  \medbreak\noindent The fact that maximal 
complexity is intimately linked to absence of triangles 
(and small cycles, more generally) is interesting and 
useful.  Theorem 1 and its corollaries point strongly to 
this fact.  We take a few lines here to intuitively argue 
why small cycles are not favored by graphs of maximal 
complexity.  Start with any tree $T$.  We are thinking of 
tossing in an extra edge $\epsilon$ in $T$ (this means keeping the 
vertex set the same and just drawing an extra edge) and 
examine what happens to the resulting complexity.  The 
resulting graph $\bar {T}=T\cup\epsilon$ has a unique cycle, call it $
C.$ It 
is immediately clear that $\bar {T}$ has exactly $c$ spanning trees, 
where $c=|C|.$ This tells us that the shorter the cycle $C$ 
is the lesser spanning trees we have!  \medbreak\noindent
Since some understanding of the structure of regular 
graphs with a minimal number of triangles exists, cf.  
[7] and [8], this reflects favorably in identifying infinite 
families of graphs of maximal complexity by way of 
Theorem 1 and its corollaries.  Large classes of 
such graphs are described in Theorem 1.6 of [8].  A more 
restricted but simple construction appears also in [7].  
We explain the details.  Let $k$ and $l$ be integers such 
that $k>l\geq 0.$ Start with a complete bipartite graph 
$K_{2k+l,2k+l}$ with vertex set $\{x_1,\ldots ,x_{2k+1}\}$ and 
$\{y_1,\ldots ,y_{2k+l}\}$.  Remove a $(l+1)-$factor from the graph 
induced by set $x_1,\ldots ,x_k,y_1,\ldots ,y_k$ and an $l$-factor from 
the graph induced by $x_{k+1},\ldots ,x_{2k+l},y_{k+1},\ldots ,y_{
2k+l}$.  Join 
$x_1,\ldots ,x_k,y_1,\ldots ,y_k$ to a new vertex $z.$ Denote by $
g(k,l)$ the 
family of graphs so obtained.  An element of $g(k,l)$ is a 
regular graph of degree $2k$ with $4k+2l+1$ vertices.  It is 
shown in [7, Theorem 2.1] that for $k\geq 2^{20}$ and 
$k\geq 2l+6\sqrt {10l}+1$ a graph in $g(k,l)$ is the sole graph with 
$4k+2l+1$ vertices of degree $2k$ having a minimal number 
of triangles (exactly $k(k-l-1)$ triangles) among all graphs 
with the same number of vertices and of the same 
degree.  When viewed in the context of Theorem 1 and 
Corollary 1 above, the results contained in Theorem 1.6 
of [8] and Theorem 2.1 of [7], provide us with infinite 
families of graphs of maximal complexity over ${\it g}_2$.  As 
has been noted, the complements of these graphs are of 
minimal complexity.  \medbreak\section*{\Large\bf4.  
Conjectures on the structure of graphs of maximal 
complexity} \medbreak A graph is said to be of {\em maximal }
{\em complexity\/} if it has a maximal number of spanning trees 
among all graphs with (at most) the same number of 
vertices and (at most) the same number of edges as 
itself.  We first state the conjectures and then present 
evidence in their support.  \medbreak\noindent {\bf Conjecture }
{\bf 1} {\em Graphs of maximal complexity, among all graphs with }
{\em at most} $n+1$ {\em vertices and at most} $m$ {\em edges, are found }
{\em iteratively as follows:  }
\vskip.3cm\noindent${\it g}_1$:  {\em They have} $n+1$ {\em vertices and} $
m$ {\em edges }
\vskip.3cm\noindent${\it g}_2:$ {\em They are as in} ${\it g}_
1$ {\em and have nearly }
{\em equal degrees }
\vskip.3cm\noindent${\it g}_3:$ {\em They are as in} ${\it g}_
2$ 
{\em and have a minimum number of triangles }
\vskip.3cm\noindent {\em For} $r\geq 4,$ \vskip.3cm\noindent${\it g}_
r:$ {\em They }
{\em are in} ${\it g}_{r-1}$ {\em and maximize} $(-1)^{r-1}Tr(C^
r)$ 
\vskip.3cm\noindent {\em After finitely }
{\em many\/} ({\em at most} $n+1$) {\em steps condition} ${\it g}_
r$ {\em involves only }
{\em cospectral graphs of maximal complexity.  }
\medbreak\noindent Settings in which regular graphs exist 
might be typically easier to handle and we reformulate 
the above conjecture for this case.  
\medbreak\noindent {\bf Conjecture 2} {\em Graphs of maximal }
{\em complexity, among all graphs with at most} $n+1$ {\em vertices }
{\em and at most} $m$ {\em edges, with} $d=\frac {2m}{n+1}$ {\em an integer, are found }
{\em iteratively as follows:  }
\vskip.3cm\noindent${\it g}_1:$ {\em They have} $n+1$ {\em vertices and} $
m$ {\em edges }
\vskip.3cm\noindent${\it g}_2:$ {\em They are as in} ${\it g}_
1$ {\em and have all }
{\em degrees equal to} $d.$ 
\vskip.3cm\noindent${\it g}_3:$ {\em They are as in }
${\it g}_2$ {\em and have a minimum number of triangles }
\vskip.3cm\noindent {\em For} $r\geq 4,$ \vskip.3cm\noindent${\it g}_
r:$ {\em They }
{\em are in} $g_{r-1}$ {\em and minimize the number of closed} $r
-${\em walks }
\vskip.3cm\noindent {\em After finitely many\/} ({\em at most} $n
+1$) 
{\em steps condition} ${\it g}_r$ {\em involves only cospectral graphs of }
{\em maximal complexity}.  \medbreak\noindent Evidently 
Conjecture 1 implies Conjecture 2.  Evidence in support 
of these conjectures is as follows.  \medbreak\noindent{\bf 1.  }
All graphs with 9 or fewer vertices and $m$ edges, for all 
possible values of $m$, satisfy Conjecture 1.  To save 
space we shall not list all these cases here.  They can 
be obtained from the authors and are also made available 
by way of an internet link.  \medbreak\noindent{\bf 2.}  An 
asymptotic heuristic, available in [6, page 157], also 
points strongly in the direction of these conjectures.  
We outline the pertinent case briefly here.  Start with a 
simple graph $G$ and ask:  what is the graph of maximal 
complexity obtained after superimposing $x$ copies of the 
complete graph on $G?$ Clearly the resulting graphs have 
multiple edges between two vertices; in fact there are $x$ 
or $x+1$ edges between any two distinct vertices.  If the 
(nonzero) eigenvalues of the Laplacian of $G$, a simple 
graph on $n+1$ vertices, are $x_1,\ldots ,x_n$, then the resulting 
graph $G+xK_{n+1}$, obtained after superimposing $x$ copies of 
$K_{n+1}$ on $G$, has eigenvalues $x_1+x,\ldots ,x_n+x.$ The logarithm 
of the complexity is \medbreak\noindent
\centerline{$\sum_iln(x_i+x)=\sum_iln(x(1+\frac {x_i}x))=c(x)+\sum_
iln(1+\frac {x_i}x)$}
\medbreak\noindent
\centerline{$=c(x)+\sum_i$$\sum_{j\geq 0}\frac {ln^{(j)}(1)}{j!}(\frac {
x_i}x)^j=c(x)+\sum_{j\geq 0}(\sum_ix_i^j)$$\frac {ln^{(j)}(1)}{j!}
(x^{-1})^j$.}
\medbreak\noindent Here $c(x)$ is a fixed function of only $x.$ 
Noting that the derivatives $ln^{(j)}(1)$ alternate in sign, this 
describes a strategy of selecting the base graph $G.$ At 
the linear level, select $G$ such that $\sum_ix_i=TrL$ is 
maximal; then, at the quadratic level, make 
$-\sum_ix_i^2=Tr(L^2)$ maximal, and so on.  This is exactly the 
strategy descibed in Conjecture 1 -- sans asymptotics.  
\medbreak\noindent Partly supportive are also the results 
in Sections 3 and 4.  These results establish the maximal 
complexity only over the class of graphs ${\it g}_2.$ In 
particular, results in Section 3 rely on early stationarity 
of chain (2) at ${\it g}_3.$ If condition ${\it g}_3$ written in Conjecture 
2 can be generally established, then large families of 
graphs of extremal complexity can be revealed by relying 
on the asymptotic structure of graphs with a minimal 
number of triangles as described in the recent paper [8].  We do 
actually believe that the same results are in fact true 
over the larger class ${\it g}_0$.  Observe that all graphs proved 
of maximal complexity over ${\it g}_2$ are as Conjecture 1 
affirms.  \medbreak\noindent
\medbreak\noindent
\centerline{APPENDIX}
\medbreak\noindent {\bf Proof of Lemma 1 \medbreak}\noindent It 
suffices to show that $f'$, the derivative of $f$ with 
respect to $x,$ is negative for $1<x<n-1.$ With $n,$ $r,$ $s$ 
(and hence $z$) fixed, we view the roots $x_1$ and $x_2$, 
written explicitly in (4), as functions of $x.$ Write 
$f=xlnx_1+(n-x)lnx_2$ and differentiate with respect to $x$ 
to obtain \medbreak\noindent
\centerline{$f'=lnx_1-lnx_2+\frac {xx_1'}{x_1}+\frac {(n-x)x_2'}{
x_2}.$}
\medbreak\noindent
We now observe that 
$xx_1'=(n-x)x_2'=\frac {-z}{2\sqrt {x(n-x)}}=\frac {-1}2(x_1-x_2)
.$ This yields 
\medbreak\noindent
\centerline{$f'=ln(\frac {x_1}{x_2})-\frac 12(x_1-x_2)(\frac 1{x_
1}+\frac 1{x_2}).$}
\medbreak\noindent Consider now the substitution 
$t=t(x)=\frac {x_1-x_2}{x_2},$ and define function 
$g(t)=ln(1+t)-\frac 12t(\frac {2+t}{1+t}),$ for $t\geq 0.$ Observe first that 
$f'(x)=g(t(x))=g(t).$ The derivative 
$g'(t)=\frac {-t^2}{2(1+t^2)}$ is negative for all $t>0.$ It follows that $
g$ 
is a strictly decreasing function of $t.$ Since also $g(0)=0,$ 
we conclude that $g(t)<0$ for all $t>0.$ This shows that 
$f'(x)=g(t(x))=g(t)<0$, for all $x,$ as asserted. 
\medbreak\noindent
{\bf Proof of Lemma 2}
\medbreak\noindent
Let $y=\sqrt {x(n-x)}$. The derivative of $y$ with respect to $x$ 
is $y'=\frac {n-2x}{2y}$. We may now write 
\medbreak\noindent
\centerline{$c(x)=y^{-3}[x(ry+(n-x)z)^3+(n-x)(ry-xz)^3]$.}
\medbreak\noindent
Upon expanding the third powers and sorting monomials 
we obtain
\medbreak\noindent
\centerline{$c(x)=nr^3+3nrz^2+2nz^3y'.$}
\medbreak\noindent
It is now easy to verify that the derivative $c'$ is 
negative. Indeed, routine calculations yield
\medbreak\noindent
\centerline{$c'(x)=2nz^3y^{\prime\prime}=\frac {-n^3z^3}{2(\sqrt {
x(n-x)})^3}<0.$}
\medbreak\noindent
This verifies the fact that function $c$ is decreasing in $x.$
\vskip1cm 
\centerline{{\bf Acknowledgement}}
\vskip.3cm\noindent\ 
We are grateful to the National Science Foundation for 
sponsoring this work under the Funding Opportunity NSF 
24-513 Emerging Mathematics in Biology.  
\vskip1cm 
\centerline{{\bf REFERENCES}}
\begin{enumerate}
\item MacDonald, I. {\em G. Symmetric functions and Hall}
{\em polynomials}, Oxford University press, 2015 

\item Brouwer, A. E. and Haemers, W. H. {\em Spectra of}
{\em graphs,\/} Springer, New York, 2012

\item van Dam, E. R. Graphs with few eigenvalues, {\em PhD}
{\em dissertation}, Tilburg University, 1996

\item Radhakrishnan, N. and Vijayakumar A. (1994), About 
triangles in a graph and its complement, {\em Discrete}
{\em mathematics,\/} 131, 205-210

\item Alon, N. (1990) The number of spanning trees in regular 
graphs, {\em Random structures and algorithms,\/} vol 1 (2), 
175-191

\item Constantine, G. M. {\em Combinatorial theory and}
{\em statistical design}, Wiley, New York, 1987

\item Lo, A. S. L. (2009) Triangles in regular graphs with density 
below one half, {\em Combinatorics, Probability and}
{\em Computing,\/} 18, 435-440

\item Liu, H., Pikhurko, O., Staden, K. (2020) The exact 
minimum of triangles in graphs with given order and 
size, {\em Forum of Mathematics, Pi,\/} Vol. 8, e8, 144 pages
doi: 10.1017/fmp.2020.7 

\item Mangin, L. et al. (2013) Neural mechanisms 
underlying breathing complexity, $PLoS$ $One$, 8(10), e75740

\item Bear, M.,Connors, B., Paradiso, M. {\em Exploring the }
{\em brain,\/} Fourth edition, Jones and Bartlett Learning, Burlington, MA, 2016

\item Higman, D. G., Sims, C. C. (1968), A simple group of 
order 44,352,000 {\em Mathematische Zeitschrift}, 105 (2), 
110-113
\end{enumerate}

\end{document}